\documentclass{article}
\usepackage{arxiv}
\usepackage[utf8]{inputenc}
\usepackage[T1]{fontenc}
\usepackage[unicode, pdftex]{hyperref}
\usepackage{url}            
\usepackage{booktabs}       
\usepackage{amsfonts}       
\usepackage{nicefrac}       
\usepackage{microtype}      
\usepackage{lipsum}		
\usepackage{graphicx}
\usepackage{amsmath}
\usepackage{amssymb}  
\usepackage{natbib}
\usepackage{doi}
\usepackage{color}

\title{The structure of Morse flows and  co-dimension one gradient flows on the sphere with holes}

	\author{
  Illia Ovtsynov\thanks{\texttt{ http://orcid.org/0009-0000-3262-2291}}
  \\
  Taras Shevchenko National University of Kyiv
	\\
	\texttt{iliaovtsynov@knu.ua}
  \And
  Alexandr Prishlyak\thanks{https://sites.google.com/view/prof-prishlyak, https://orcid.org/0000-0002-7164-807X}
  \\
  Taras Shevchenko National University of Kyiv
	\\
	\texttt{prishlyak@knu.ua}
}

\renewcommand{\shorttitle}{The structure of Morse flows and  co-dimension one gradient flows on the sphere with holes }

\newtheorem{theorem}{Theorem}

\newtheorem{corol}{Corollary}

\newtheorem{example}{Example}

\hypersetup{
pdftitle={The structure of Morse flows and  co-dimension one gradient flows on the sphere with holes},
pdfauthor={Alexandr O. Prishlyak},
pdfkeywords={bifurcation, topological equivalence, structural stability},
}

\begin{document}
\maketitle

\begin{abstract}
We describe all possible topological structures of typical one-parameter bifurcations of gradient flows  on the 2-sphere with holes in the case that the number of singular point of flows is at most six. To describe structures, we separatrix diagrams of flows. The saddle-node singularity is specified by selecting a separatrix in the diagram of the flow befor the bifurcation and the saddle connection is specified by a separatrix, which conect two saddles.

\end{abstract}


\section*{Introduction}
\label{sec:1}

In 1970, Palis and Smale \cite{Palis1970} established that Morse--Smale gradient vector fields exhibit structural stability within the collection of all gradient fields $\operatorname{Grad}(M)$ with a $C^k$ topology ($k \ge 3$) on a smooth manifold $M$. Additionally, the set of Morse--Smale vector fields $\operatorname{MS}(M)$ is both open and dense within $\operatorname{Grad}(M)$. Gradient vector fields of codimension 1 are those that occur in standard 1-parameter families of gradient vector fields and do not qualify as Morse--Smale fields. In 1983, Palis and Takens \cite{Palis1983} demonstrated that these fields can be categorized into two types: SN-fields, which feature a single saddle-node singularity, and SC-fields, which include a saddle connection. Aside from these singularities, they possess the same characteristics as Morse--Smale gradient fields. Specifically, they maintain structural stability within the set $\operatorname{Grad}(M) \setminus \operatorname{MS}(M)$, and they also constitute an open and dense subset in this space.

On closed manifolds, a vector field consistently produces a flow. In the case of manifolds with boundaries, a vector field defines a flow only if it is tangent to the boundary. A flow is considered gradient, Morse--Smale, or of codimension one if the generating vector fields possess these respective properties. The description of Morse--Smale flows on manifolds with boundary can be found in \cite{lp90}, while flows of codimension one are addressed in \cite{Kibkalo2022}. The structural stability of these flows within the relevant spaces, as well as the openness and overall density of their sets in these spaces, has also been established.

Complete topological invariants can be established for the topological classification of flows. These invariants often take the form of graphs that possess specific properties and are augmented with additional information. Furthermore, two flows are considered topologically equivalent if there exists an isomorphism between their graphs that maintains this information. The effectiveness of these invariants is determined by the ability to identify all structures of dynamical systems given a specific number of singular points.

For Morse--Smale flows, a separatrix diagram serves as a complete topological invariant \cite{Fleitas1975, Leontovich1955, Oshemkov1998, Peixoto1973}.

We can refine this invariant for Morse flows (which are Morse--Smale flows lacking closed trajectories) that are topologically equivalent to the gradient flows of Morse functions. This invariant manifests as a graph embedded in a surface, where the vertices represent sources and the edges correspond to the stable manifolds of saddles.

A rotation system \cite{GT87} is frequently utilized to define a graph embedding within a surface. This serves as a Peixoto invariant \cite{Peixoto1973} for Morse flows. When dealing with a flow that has a single source, chord diagrams are a convenient means of representation \cite{Kybalko2018}. To create a chord diagram, codes are employed, which consist of the numbers representing the chords whose endpoints connect while moving around the circle. However, the potential for ambiguity in numbering the chords or selecting a starting point can make it challenging to compare the resulting invariants.

Additionally, complete topological invariants for flows exhibiting gradient dynamics have been developed in two dimensions \cite{prishlyak2020three} and three dimensions \cite{Prishlyak2002beh2}.
In  \cite{Prishlyak2020, Prishlyak2021}, flows with collective dynamics on the sphere were studied, which, in addition to Hamiltonian regions, are Morse flows on the sphere with holes. In this case, the saddle points on the boundary of the gradient region contain three hyperbolic sectors, unlike the flows we are considering, which contain two hyperbolic sectors for such points.

The topological structure of Morse flows can also be studied using Lyapunov functions, which are Morse functions. Their classification was obtained in \cite{hladysh2019simple, Hladysh2019, Khohliyk2020, Kravchenko_2019, prish2001,  prish2015, Prishlyak2000, Reeb1946, Hatamian2020}. 



The purpose of this paper is constructing complete invariant for Morse and gradient codimension one flow on the sphere with holes which resembles a chord diagram for Morse flows. This invariant has a marked point in the diagram, which allows us to define clearly a number code of the flow.
The invariants we constructed (the distinguishing graph and the flow code ) are generalizations of the distinguishing graph of Peixoto and the Oshmankov-Shark code, which were developed for Morse flows on closed surfaces.

The results of this article were represented in \cite{ShV2025, AGMA2025, Prish-Ovt_2025}. 







\section{Topological Invariants of Morse Flows}

Typical vector fields on compact 2-manifolds are Morse-Smale fields. Among gradient fields, these are Morse fields or gradient-like Morse-Smale fields that do not contain closed trajectories. They satisfy three properties:

1) The critical points are non-degenerate;

2) There are no separatrix connections belonging to the interior of the 2-manifold;

3) The $\alpha$-limit ($\omega$-limit) set of each trajectory is a critical point.

Recall that the neighborhood of an isolated critical point, apart from the center, can be divided by trajectories into sectors, which can be of one of three types: 1) elliptic (the $\alpha$- and $\omega$-limit sets of the trajectories are this critical point), 2) hyperbolic (neither the $\alpha$- nor the $\omega$-limit sets of the trajectories are this critical point), 3) parabolic (only one of the $\alpha$- or $\omega$-limit sets is this critical point). By \textit{separatrix}, we mean an internal trajectory that intersects the boundary of a hyperbolic sector of some critical point.

The \textit{separatrix diagram} of a gradient flow (vector field) is a surface on which critical points, boundary trajectories, separatrices, and a specified orientation of these trajectories are highlighted. For convenience, we   depict separatrices that start at a saddle and end at a sink (\textit{unstable separatrices}) in green, and those that start at a source and end at a saddle (\textit{stable separatrices}) in red.

The separatrix diagram is a complete topological invariant of a Morse flow on a compact surface. However, this invariant can be simplified by removing internal sinks and unstable (green) separatrices. The remaining critical points are 0-cells, the rest of the boundary trajectories and red (stable) separatrices are 1-cells, and the components of the complement to the 0- and 1-cells are 2-cells. Each 2-cell has a unique sink: either one of the critical points on the boundary or an interior point. Since, up to homeomorphism, saddle points on the boundary of a cell can be connected by curves without intersection with the sink uniquely, the cell decomposition we described carries the same information about the flow structure as the separatrix diagram.

\subsection{Distinguishing Graph of the Morse Flow}

As with Morse flows on closed manifolds, the separatrix diagram is a complete topological invariant of such flows. It consists of an oriented graph embedded in the surface. The vertices of the graph are the critical points, and the edges are the separatrices and trajectories belonging to the boundary. The orientations of the edges are determined by the direction of movement along the corresponding trajectories. The separatrices entering the saddles are called stable, while those exiting them are called unstable. To define the flow, the diagram can be simplified by removing unstable separatrices from the graph and collapsing internal saddles (replacing an internal saddle and its two stable separatrices with a single edge), as well as collapsing b-saddles on the boundary. The remaining critical points on the boundary are colored in three colors: 1) sinks in white, 2) a-saddles in black, and 3) sources in gray.

For the obtained graph embedded in a surface, we can use a rotation system. To define the cyclic order of half-edges at each source vertex, we   conduct an additional construction. For each source, we   take a sufficiently small regular neighborhood with a smooth boundary, which each trajectory intersects at most at one point. We   color the boundary of the neighborhood in a second color (dashed). We   remove the interiors of these neighborhoods, and the intersections of their boundaries with the separatrices and the trajectories of the boundary   form new vertices of the graph, which we   color black (the second color).

Edges that intersect with the separatrices are colored in the first color (thin black), while edges belonging to the boundary are colored in the third color (thick black). Thus, we obtain a graph where all vertices are colored in two colors and edges in three colors. We   refer to this graph as the distinguishing graph of Morse flow.

A \textit{boundary cycle} of the distinguishing graph are  called a component of the connectivity of the graph obtained by removing edges of the first color. The orientation of the surface defines the orientation of each boundary cycle.

\begin{figure}[ht!]
\center{\includegraphics[width=0.40\linewidth]{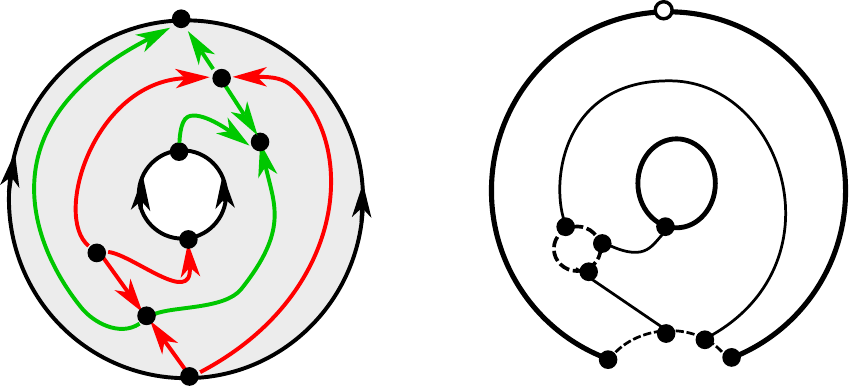}
}
\put (-33,13) {\small 2}
\put (-40,1) {\small 1}
\put (-58,15) {\small 1}
\put (-67,35) {\small 2}
\put (-38,36) {\small 3}
\put (-53,33) {\small 3}
\caption{Separatrix diagram and destiguished graph for Morse flow with code $\{12]0[\}[123]\{3\}.$}
\label{dg}
\end{figure}
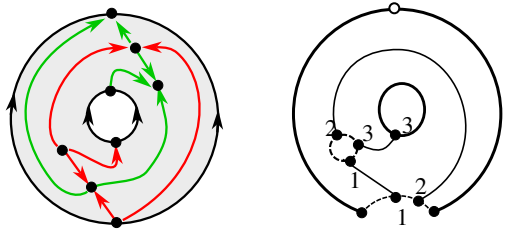

An example of a distinguishing graph is depicted in Fig. \ref{dg}. In this and other figures, the orientation of the outer boundary cycle is defined by counterclockwise movement, while the orientations of the inner boundary cycles are defined by clockwise movement.

\begin{theorem} \label{t1} Two Morse flows are topologically equivalent if and only if there exists an isomorphism between their distinguishing graphs that preserves the colors of vertices and edges, and either preserves the orientations of all boundary cycles or simultaneously reverses the orientations of all boundary cycles.
\end{theorem}
\textit{Proof.}
Necessity. A topological equivalence $h$, which is a homeomorphism, induces a bijection between special points, separatrices, and boundary components, thus defining an isomorphism of graphs. By construction, the colors of the corresponding vertices and edges are the same. If $h$ preserves the orientation of the surface, then the orientations of all boundary cycles are preserved; if $h$ reverses the orientation of the surface, then the orientations of all boundary components are simultaneously reversed.

Sufficiency. We double the distinguishing graph $G$. To do this, we identify the points of edges of the third color on the disconnected union $G \sqcup G$. In this way, we construct the distinguishing graph of the doubled flow on a closed surface. The resulting flows on the double are $\Omega$-stable flows, and the equivalence of distinguishing graphs (an isomorphism that preserves colors and orientations) defines the topological equivalence of the flows. 
The restriction of these homeomorphisms to one half of the surface is the sought-after topological equivalence.

\subsection{Morse Flow Code}

In this section, we   consider Morse flows with a distinguished vertex on a distinguishing graph and a given orientation of the surface (which is equivalent to specifying an oriented orientation on one (and thus, on all) of the boundary cycles $C_1, C_2, \ldots$. All saddle points will be numbered. Their stable manifolds   receive the same numbers. The code   consist of lists of the numbers of separatrices encountered during traversal according to the orientation of the boundary components of the surface (the list is enclosed in curly braces) or the boundaries of the neighborhoods of sources (the list in square brackets).

We start the code with an opening curly brace if the vertex belongs to the boundary or the neighborhood of a source on the boundary, or with an opening square bracket if the vertex belongs to the boundary of an internal saddle. If the distinguished vertex belongs to a separatrix, we number the corresponding saddle (and its separatrices) as 1 and include it in the code. If at the vertex, while moving according to the orientation of the boundary cycle $C_1$, there is a color change from the third to the second, it corresponds in the code to the symbol $[$, and if from the second to the third, then the symbol $]$. We apply this rule for all vertices of this type. We move from the distinguished vertex to the next one. If it belongs to the separatrix of a saddle that has not yet been numbered, we assign it the next number (number 1 if no saddle has been numbered yet) and include it in the code. If this vertex corresponds to an already numbered saddle or the symbols $[,]$, we include that number or symbol in the code. We continue this process until we reach the distinguished point. We close the list with the appropriate curly or square bracket.

If a distinguishing graph has one boundary cycle, the process of code construction is complete. Let's assume there are several such cycles. Since the distinguishing graph is connected, there exists an edge of the first color that connects the first boundary cycle to others. Among all such edges, we   select the edge $e$ with the smallest number. The corresponding boundary cycle is the second cycle $C_2$, and the end of $e$ that belongs to $C_2$ is the distinguished vertex on $C_2$. For cycle $C_2$, the process of listing is the same as for $C_1$. Following the same principle, we move on to the next boundary cycle. The code writing is completed when we exhaust all boundary cycles. If we encounter a white vertex while traversing the boundary cycle, we   place a 0 at the corresponding position. If two numbers appear consecutively in the list, we   separate them with a comma. If there are fewer than 10 seats in the stream (all numbers are single-digit), we will not use commas in the code.

Thus, we have two types of lists in the code: 1) internal lists, which are located in square brackets and are not contained within curly brackets (and do not contain other brackets), and 2) boundary lists, which are located in curly brackets (and do not contain other curly brackets). If the boundary list contains square brackets, we   cyclically rearrange its elements so that it starts with an opening square bracket. Sub-lists of the code are defined as sequences in the modified boundary lists that are placed between the opening and closing square brackets, provided that they do not contain other square brackets.

\begin{example} Let a vertex corresponding to the lower separator of the lower saddle be highlighted on the distinguishing graph in Fig. \ref{dg}, with the orientation of the outer boundary cycle counterclockwise. Then the code  looks like: $\{12]0[\}[123]\{3\}.$ It consists of three lists and contains one sublist.
\end{example}

\begin{theorem}  \label{t2} There exists an isomorphism between distinguishing graphs that maps the highlighted vertex to the highlighted one and preserves the colors of the vertices, the colors of the edges, and the orientations of the boundary cycles if and only if they have the same codes.
\end{theorem}
\textit{Proof.} Necessity. The choice of the highlighted vertex and the orientation of the boundary cycles allows for a unique numbering of the vertices of the distinguishing graph corresponding to the saddle special points of the flow. The correspondence of the colors of the edges of the distinguishing graph defines a correspondence between square and curly brackets in the code. Therefore, the codes of the graphs are the same.

Sufficiency. Let the flows have the same codes. We   construct an isomorphism between their distinguishing graphs. The boundary cycles of the graphs correspond to the lists of the code. We   define homeomorphisms between the cycles that preserve orientations and map vertices with the same numbers and in the same positions in the code to vertices with the same numbers and positions. We   extend the constructed homeomorphisms to the internal edges. This can be done because, by construction, equality of numbers implies consistency of the constructed homeomorphisms with the ends of the internal edges, and the mapping of the ends of the segments extends to their homeomorphisms. Thus, we have an isomorphism of distinguishing graphs. Since the second and third colors of the edges are determined by the placement of square brackets, the constructed isomorphism of the graphs preserves the colors of the edges. The numbered 0 vertices have a different color than those numbered with positive integers. Therefore, the colors of the vertices are also preserved.

We   refer to a \textit{source list} as a list of numbers that are placed between opening and closing square brackets or between such brackets in some cyclic permutation of symbols within curly brackets.

For the code, we   write a system of cycles. To form a cycle, we   fix a certain number (let's denote it as $a_1$) in one of the lists: this is the first number in the cycle. The next number in the cycle is (cyclically) the next number in the fixed list, which we   call $a_2$. Next, we find another list in the code that contains the number $a_2$ (or the same list with $a_2$ in a different location) and move to it. We take (cyclically) the next number $a_3$ and repeat the procedure. The formation of the cycle is completed when we return to the initial $a_1$ at the same position in the same list from another $a_1$. If a 0 is encountered while listing, we include it in the cycle and continue taking the next number from the same list. The cycle corresponds to the boundary of 2 cells in the cellular decomposition.

\begin{theorem}
The Morse code has the following properties:

1. Each number, except for 0, appears in the code twice and is contained in at least one of the source lists.

2. Each cycle contains at most one 0.

3. For any two lists in the code, there exists a sequence of lists such that each pair of adjacent lists in this sequence contains the same numbers.

4. If $V$ is the number of lists (pairs of curly brackets and pairs of square brackets that do not lie within curly brackets), $E$ is the number of numbers (excluding 0), and $F$ is the number of cycles, then the Euler formula holds:

$$V-E+F=2.$$

\end{theorem}
Proof.
1. Non-zero numbers of the code correspond to stable manifolds of saddle points that lie in the interior of the surface. There are two possibilities: 1) an internal saddle, 2) a saddle on the boundary. In the first case, the stable manifold consists of two separatrices, each of which starts at a source. We have two identical numbers in the lists of sources. In the second case, the stable manifold is a separatrix. One of the corresponding numbers is in the list of sources, while the other is the number of the saddle on the boundary.

2. Since a cycle corresponds to a cell (face) of the distinguishing graph, and each face contains only one sink, we have either one 0 in the cycle if this sink belongs to the boundary, or no 0 if the sink is an interior point.

3. This condition is equivalent to the connectivity of the distinguishing graph, which follows from the connectivity of the surface.

4. Using Morse flow, we   construct a cell decomposition of the sphere. For this, we collapse each boundary component to a point. The 0-cells are the sources of the flow and the boundary components collapsed to a point. They correspond to lists. The 1-cells are the stable manifolds of saddle points. They correspond to numbers. The 2-cells correspond to cycles. Then for the Euler characteristic of the sphere, $\chi(S^2) = V - E + F = 2.$

\begin{theorem}
Let a set of lists be given, consisting of integers from 0 to n>0, where each list is written either 1) in curly braces or 2) in square brackets; lists in curly braces may also contain sublists written in square brackets (with subsequent cyclic permutations of the elements of the list). If these sets of lists satisfy properties 1---4 of the previous theorem, then there exists a Morse flow on a compact surface of genus 0 (a sphere with holes) for which they serve as the code of the flow.
\end{theorem}
\textit{Proof.}
We   construct a distinguishing graph from the code. For this, we take the disjoint union of $n$ circles $S^1$, where $n$ is the number of lists. Thus, each list corresponds to a circle. We   assign an orientation (direction of movement) on each circle arbitrarily. On each circle, we   select points (vertices) that correspond to the numbers and the inner square brackets, so that when traversing according to the orientation, their order matches the order of the corresponding elements of the list. We   connect vertices with the same numbers by segments (without pointwise intersection). We   color the edges of these segments in the first color, the edges of the circles corresponding to lists and sublists in curly braces in the second color, and the remaining edges in the first color. The vertices corresponding to 0 is colored in the second color, while the rest of the vertices are colored in the first color. After doubling the distinguishing graph, as in Theorem 1, we   construct a flow on the doubled surface. 
The restriction of this flow to half of the surface is the desired Morse flow.

By construction, the code of the flow depends on the initial special point and the orientation of the surface (boundary). Changing the initial point alters the numbering of the vertices, meaning a substitution occurs in the set of vertex numbers. Additionally, in each list, changing the initial point may lead to a cyclic permutation of the elements of the list. The order of the lists in the code may also change. If the orientation of the surface changes to the opposite, this  results in the order of the elements of each list being reversed (with possible cyclic permutations).

\section{Typical One-Parameter Bifurcations of Gradient Flows on Spheres with Holes}

In typical one-parameter families of fields, one of the three conditions of Morse flows may be violated. Violation of the first condition leads to a bifurcation at a singular point of codimension 1, while violation of the second condition results in the emergence of a saddle connection. The third condition cannot be violated, as per the Poincar\'e-Bendixson theorem, the $\alpha$-limit set ($\omega$-limit set) of each trajectory on the sphere is either a singular point, a cycle, or a polycycle. Since gradient fields do not have cycles or polycycles, this set is a singular point.

All possible typical bifurcations in one-parameter families of gradient flows are described in \cite[Theorem I]{Kibkalo2022}. These bifurcations include: 1) saddle-node (internal) SN, 2) boundary saddle-node BSN, 3) degenerate saddle on the boundary (half-saddle) HS, 4) degenerate node on the boundary (half-node) HN, 5) double saddle on the boundary BDS, 6) saddle connection (internal) SC, 7) half-boundary saddle connection (semi-boundary) HSC, 8) boundary saddle connection BSC. Bifurcations 1)--5) are often referred to as local bifurcations, while 6)--8) are considered global.

We   examine the local bifurcation in the vicinity of $U$ of the singular point ($x=0$, $y=0$) as a family of vector fields depending on the parameter $a$: $V=V(x,y,a)$, $a \in [-1,1]$, with $a=0$ being the bifurcation point. We   choose the direction of the parameter $a$ such that for $a<0$, the vector field in the vicinity of $U$ has two singular points. In all cases, the vector field $V(x, y, 0)$ can be obtained from $V(x, y, -1)$ by collapsing to a trajectory that connects the singular points. Thus, at $a=0$, we have a flow (vector field) of codimension 1.

To define a flow with a special point of codimension 1, one can either A) identify a trajectory in the topological invariant of the vector field before bifurcation that converges to a point, or B) indicate the location of the bifurcation, the type of bifurcation, and which other elements of the separatrix diagram changed during the bifurcation on the flow after bifurcation. The first method is simpler but requires defining the flow with one more special point than the flow of codimension 1. Thus, to find all structures of flows of codimension 1 with $n$ special points using the first method, one must know the structures of all Morse flows with $n+1$ special points. For the second method, it is sufficient to know the structures of Morse flows with $n$ and $n-1$ special points. Therefore, for each type of flow 1)--5), we   construct two types of codes: A-code and B-code.

\subsection{Saddle-node (internal) SN}

\begin{figure}[ht!]
\center{\includegraphics[width=0.75\linewidth]{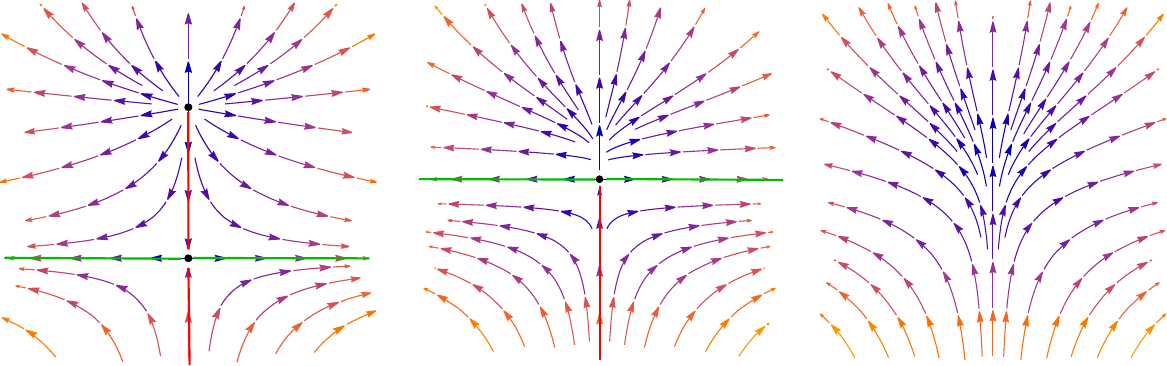}
}
\caption{saddle-nod bifurcation $SN_{+}$
}
\label{bifsn}
\end{figure}

The saddle-node bifurcation ($SN$) occurs when the gradient bifurcation has an internal degenerate special point of codimension 1, and one of the eigenvectors of the derivative matrix is equal to 0 while the other is non-zero. It can be described by the equation (normal form) $V(x,y,a)=\{x,y^2+a\}$ if the node is a source (bifurcation $SN_{+}$), and by the equation $V(x,y,a)=\{-x,-y^2+a\}$ if the node is a sink (bifurcation $SN_{-}$). Here, $a$ is a parameter. The bifurcation $SN_{+}$ is depicted in Fig. \ref{bifsn}. If the direction of motion along all trajectories is reversed in the bifurcation $SN_{+}$, we obtain $SN_{-}$.

For $a<0$, we obtain the flow before the bifurcation, for $a>0$ we obtain the flow after the bifurcation, and for $a=0$ we obtain the flow at the moment of bifurcation (the flow of corank 1). It should be noted that the saddle-node bifurcation corresponds to a typical one-parameter catastrophe (deformation) of the functions $f(x,y,a)=\pm \frac{1}{2} x^2+\frac{1}{3}y^3+ay$:
$$ \text{grad} \ (\pm \frac{1}{2} x^2+\frac{1}{3} y^3+ay)= \{\pm x, y^2+a\}.$$

To define the bifurcation $SN_{+}$, we can use either 1) the Morse flow before the bifurcation ($a=-1$) or 2) the Morse flow after the bifurcation ($a=-1$).

1) Consider the stable manifold (in Fig. \ref{bifsn}, the red separatrix) for the special saddle-source point. It cannot start at this same special point since the gradient field does not contain loops. This means that the stable saddle manifold involved in the saddle-source bifurcation does not form a loop. Therefore, to determine the saddle-source bifurcation on the distinguishing graph of the Morse flow before the bifurcation (at $a=-1$), we need to highlight a vertex on the boundary cycle of the source. The orientation of the surface can be defined in two ways that are not related to the bifurcation. Thus, there are two A-codes for $SN$.

For example, for the flow in Fig. \ref{dg} for the first separatrix emanating from the inner node, we have the following A-code: $SN_+^A[123]\{12]0[\}\{3\}.$

2) Let $\gamma$ be a separatrix that belongs to $SN_{+}$, $\alpha (\gamma)=z$. The separatrices of saddle points that emerge from a saddle-node bifurcation start at $z$. To represent the bifurcation of $SN_{+}$ on the separatrix diagram of the Morse flow, these trajectories must be highlighted after the bifurcation. In the code, we   highlight them using parentheses.

\begin{example} For the same bifurcation described for the A-code, the B-code looks like: $SN_+^B\{(12)1]0[\}\{2\}.$
\end{example}

At the moment of bifurcation ($a=0$), the separatrix diagram of the flow defines the structure of the bifurcation.

Similarly, a bifurcation of $SN_{-}$ is possible if the unstable manifold of the saddle before the bifurcation does not form a loop. In the separatrix diagram of the Morse flow before the bifurcation, the unstable separatrix should be highlighted. In the cellular decomposition, the corresponding 1-cell is located between two different 2-cells. To define such a bifurcation, a pair must be highlighted: one cell whose closure belongs to the interior of the surface, and one of the two adjacent 2-cells. Moreover, the adjacent 2-cell should not contain a sink on its boundary and should only be adjacent to the 1-cell on one side.

Additionally, to define the structure of such flows, one can reverse the directions of movement along all trajectories of the flow and use the constructed codes for $SN_{+}$-flows.

For example, for the separatrix that emerges from the first saddle in Fig. \ref{dg}, the internal sink codes will be as follows: $SN_-^A[123]\{1]0[3\}\{2\}$ and $SN_-^B\{(12)]0[2\}\{1\}$.

\textbf{Note 1.} The A-code of the SN-flow begins with a square bracket, while the separatrix that is being pulled corresponds to 1. In the B-code, the round brackets are either in the list of internal sources, that is, in the list within square brackets, or in one of the sublists, and they appear at the beginning of the code after a square or curly bracket.

\textbf{Note 2.} The bifurcations of saddle-node in Morse theory correspond to the operation of collapsing points or collapsing complementary handles of adjacent indices.

\subsection{Boundary Saddle-Node BSN}
A boundary saddle-node bifurcation occurs when a singular point lies on the boundary, the eigenvector belonging to the boundary is zero, and the transverse to the boundary is non-zero. Its normal form is given by $V(x,y,a)=\{x^2+a, y\}$, $y\ge 0$, if the node is a source (bifurcation $BSN_{+}$), and by the equation $V(x,y,a)=\{-x^2+a, -y \}$, $y\ge 0$, if the node is a sink (bifurcation $BSN_{-}$).

Note that the boundary saddle-node bifurcation corresponds to the deformation of the functions $f(x,y,a)=\pm \frac{1}{3} x^3+ax \pm \frac{1}{2}y^2$.

\begin{figure}[ht!]
\center{\includegraphics[width=0.75\linewidth ]{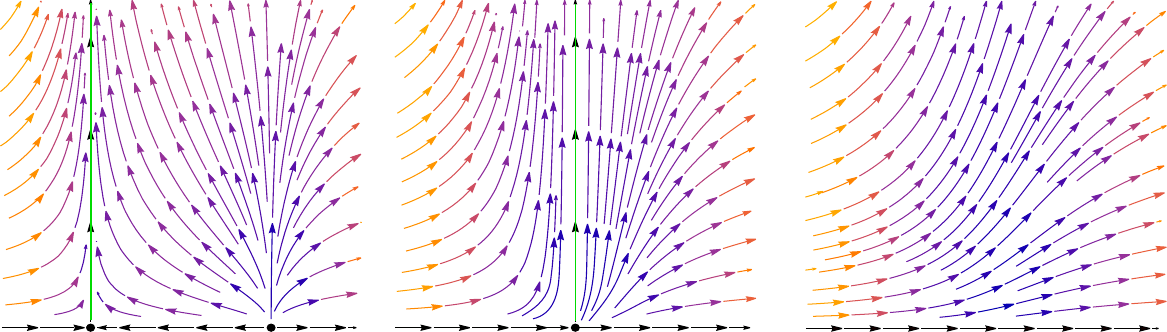}}
\caption{  $BSN_{+}$ bifurcation }
\label{s-d}
\end{figure}

The bifurcation $SN$ is $Z_2$ invariant (has an axis of symmetry), and $BSN$ can be viewed as a factorization of $SN$ under the action of this group. In other words, if at any moment we cut the phase diagram of the flow $SN$ ($V(x,y,a)=\{pm x,\pm y^2+a\}$) along the axis of symmetry ($x=0$), then bifurcation $BSN$ occurs in each part.

For $BSN_{+}$, the source and saddle merge into a point.
In Fig. \ref{s-d} a) the flow before the bifurcation is shown ($a=-1$), in Fig. \ref{s-d} b) -- the flow at the moment of bifurcation ($a=0$), and in Fig. \ref{s-d} c) -- the flow after the bifurcation ($a=1$).

For $BSN_{-}$, the saddle and sink merge into a point, which then disappears. 


As with the internal saddle-node, the structure of the bifurcation $BSN$ can be defined by a separatrix diagram at the moment of bifurcation.

To write the A-code, we   define the bifurcation using a distinguishing flow graph before the bifurcation, where the boundary trajectory that is attracting is highlighted. The highlighted point is the saddle, and the orientation of the boundary is determined by the direction of movement from the highlighted saddle to the node.

It should be noted that the boundary component containing the attracting trajectory must include at least 4 special points (if there are fewer points, then after the bifurcation we   obtain a closed trajectory or loop on the boundary, which is not possible for gradient bifurcations).

\begin{example} In Fig. 11 ??? (9), we highlight the boundary trajectory with the lower source and the left saddle. This trajectory corresponds to the first left square bracket ([) in the code. Therefore, the code looks like this:
$BSN_+^A\{[1]0[]\}[1]$.
\end{example}

For the B-code on the distinguishing Morse flow graph after the bifurcation, we need to highlight the boundary edge $a$. If $\alpha(a)$ is a source, then we need to highlight the separatrices that emerged from the BSN at the moment of bifurcation. The highlighted separatrices in the code are placed immediately after the bracket [, corresponding to edge $a$. If $\alpha(a)$ is a saddle and $\gamma$ is the separatrix of this saddle, then $\alpha(\gamma)$ is a source, and we need to highlight the separatrices of this source that lie on one side of $\gamma$, which is also edge $a$. In this case, $\gamma$ is the first highlighted separatrix.

In the B-code, just like for $SN$, we highlight the numbers (in parentheses) of the separatrices that emerge from the attracting source or the beginning $\gamma$. If there are no such separatrices, we leave an empty space in parentheses.

\begin{example} For the same bifurcation described in the previous example, the B-code is: $BSN_+^B\{[(1)]0\}[1]$.
\end{example}

\begin{example} For the Morse flow 6) in Fig. 13 ???, we highlight the right external semicircle. If the right red separatrix is highlighted, the B-code looks like this: $BSN_+^B\{10\}\{12]0[(3)\}$. If no separatrices are highlighted, the B-code is: $BSN_+^B\{10\}\{12]0[3( ) \}$.
\end{example}

To construct the codes $BSN_{-}$, as well as for $SN$, we use the reversed flows.

The next two cases of the special point of corank 1 on the boundary correspond to the situation where the eigenvector of the linear part matrix belonging to the boundary is non-zero, while the transversal to the boundary is zero. Then the normal form of the vector field can be written as $V(x,y,a)=\{\pm x, \pm y^2+ay \}$, $y\ge 0$. These bifurcations can be described as those where the separatrix, which collapses to a point, has one end on the boundary and the other is an interior point of the surface.

\subsection{Degenerate saddle on the boundary HS}

\begin{figure}[ht!]
\center{\includegraphics[width=0.75\linewidth]{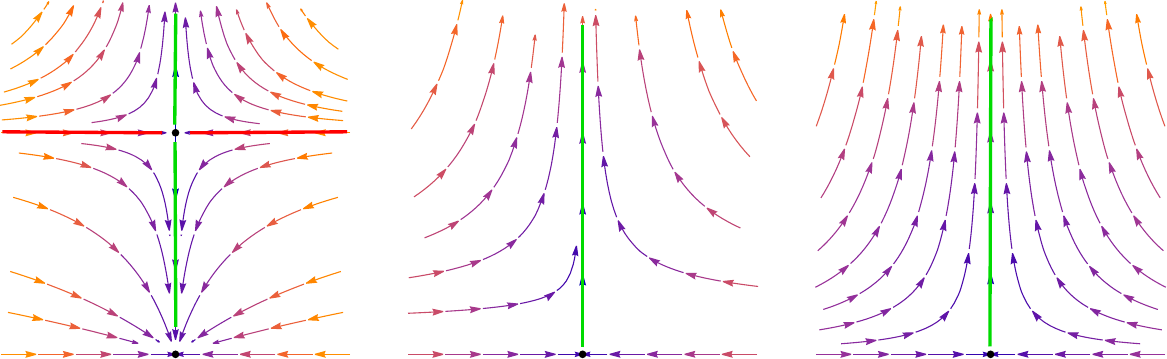}}
\caption{$HS_-$ bifurcation}
\label{HS}
\end{figure}

There are two types of degenerate saddles on the boundary that arise at $a=0$ for such bifurcations: $HS_-$ with the normal form $V(x,y,a)=\{- x, y^2+ay \}$, $y\ge 0$ (see Fig. \ref{HS}) and $HS_+$ with the normal form $V(x,y,a)=\{ x, -y^2 -ay \}$, $y\ge 0$ and $HS_+$. They are topologically equivalent to the bifurcations $V(x,y,a)=\{- x, y^3+ay \}$ and $V(x,y,a)=\{ x, -y^3 -ay \}$, $y\ge 0$, which are symmetric with respect to the axis $y=0$. $HS_+$ is the reverse of $HS_-$.

Note that the bifurcation of the boundary saddle-node corresponds to the deformation of the functions $f(x,y,a)=\mp \frac{1}{3} x^3+ax \pm \frac{1}{2}y^2$.

To construct the A-code of the $HS_+$ flow, one must highlight the source on the boundary and the separatrix that emerges from it on the flow diagram before the bifurcation.

For example, for the flow in Fig. \ref{dg} and the first separatrix, the code is recorded as $HS_+^A\{12]0[\}[123]\{3\}.$

To construct the B-code of the $HS_+$ flow, one must highlight one of the saddles on the flow diagram after the bifurcation, which was degenerate at the moment of bifurcation.

For instance, if for the flow in Fig. \ref{dg} saddle 3 on the boundary is highlighted, then after renumbering the saddles, the B-code takes the form: $HS_+^B\{1\}[123]\{23]0[\}.$

The A- and B-codes of the $HS_-$ flows are constructed as the corresponding codes of the reversed $HS_+$ flows.

\subsection{Degenerate Node on the Boundary HN}

\begin{figure}[ht!]
\center{\includegraphics[width=0.75\linewidth ]{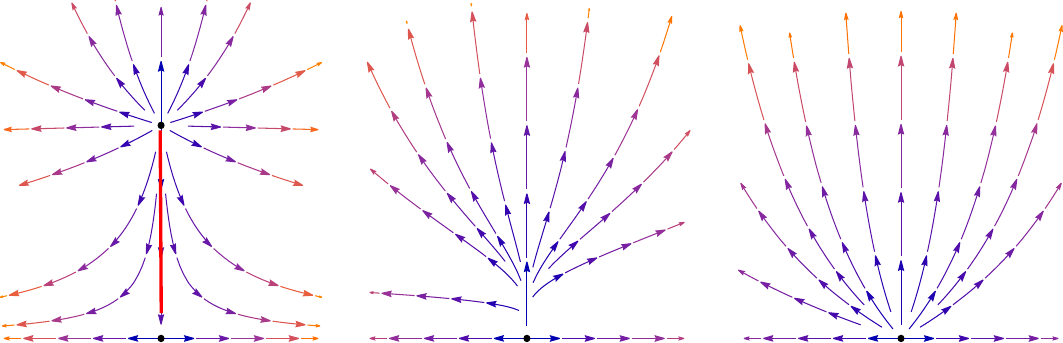}}
\caption{ HN bifurcation }
\label{HN}
\end{figure}

There are two types of degenerate nodes on the boundary that arise at $a=0$ for such bifurcations: the source $HN_+$ with the normal form $V(x,y,a)=\{ x, y^2 -ay \}$, $y\ge 0$ (see Fig. \ref{HN}) and the sink $HN_-$ with the normal form $V(x,y,a)=\{- x, - y^2+ay \}$, $y\ge 0$. Topologically equivalent to them are the bifurcations $V(x,y,a)=\{x, y^3 - ay \}$ and $V(x,y,a)=\{ -x, -y^3 + ay \}$, $y\ge 0$, which are symmetric with respect to the axis $y=0$. $HN_-$ is the reverse of $HN_+$.

It should be noted that the bifurcation of the degenerate node on the boundary corresponds to the deformation of the functions $f(x,y,a)=\pm \frac{1}{2} x^2 \pm \frac{1}{3}y^3\mp \frac{1}{2} ay^2 $.

To construct the A-code of the $HN_+$ flow, one must highlight the saddle on the boundary and the separatrix entering it on the flow diagram before the bifurcation.

For example, for the flow in Fig. \ref{dg} and the third separatrix after renumbering the saddles, the A-code takes the following form: $HN_+^A\{1\}[123]\{23\}0[\}.$

To construct the B-code of the $HN_+$ flow, one must highlight one of the sources on the boundary after the bifurcation, which was degenerate at the moment of bifurcation.

For instance, if the lower source on the boundary is highlighted for the flow in Fig. \ref{dg}, the B-code is recorded as: $HN_+^B\{[12]0\}[123]\{3\}.$

The A- and B-codes of the $HN_-$ flows are constructed as the corresponding codes of the reversed $HN_+$ flows.

\subsection{Double saddle on the boundary BDS}

\begin{figure}[ht!]
\center{\includegraphics[width=0.75\linewidth ]{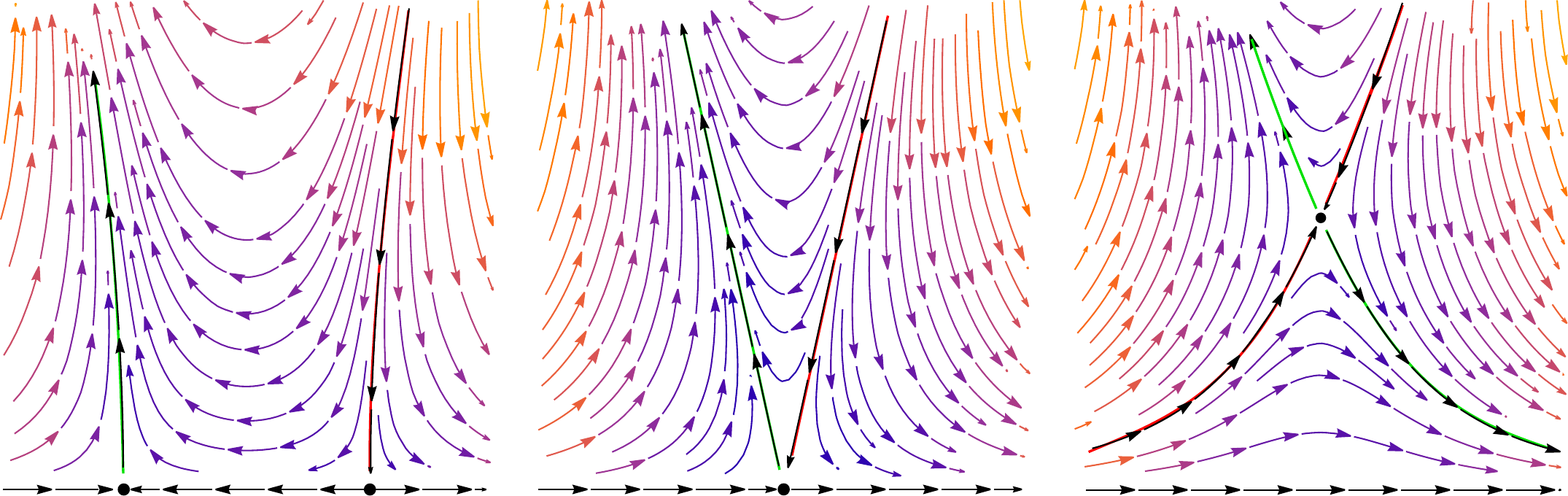}}
\caption{ BDS bifurcation}
\label{s-s}
\end{figure}

In this bifurcation, both eigenvalues of the linear part matrix are equal to 0. It occurs when the trajectory of the Morse field belonging to the boundary and connecting two saddles collapses to a point. Its normal form is given by
$V(x,y,a)=\{ x^2-y^2+a, -2xy \}$ (see Fig. \ref{s-s}). It can be realized as the gradient of such a family of functions: $f(x,y,a)=\operatorname{Re} (x+ \operatorname{i} y)^3+ax$..

At this bifurcation, an internal special point appears, since when doubling the vector field on the sphere, according to the Poincar\'e-Hopf theorem, the sum of the Poincar\'e indices does not change during bifurcation, so two saddles on the boundary (index -1) produce one saddle in each of the half-planes.

To define the A-code among the merging saddles, it is necessary to highlight the one that includes the separatrix, and the direction of the boundary containing it should be set as the direction of movement along the trajectory between the saddles.

For example, let's say that in diagram 6) in Fig. 11, the trajectory between the lower and left saddle is highlighted. Then the A-code looks like this: $BDS^A\{120\}[1][2]$.

To write the B-code for the flow after bifurcation, it is necessary to highlight the saddle (separatrix) and the trajectory on the boundary to which it adheres. In this case, the saddle and the boundary trajectory must belong to the same cell (one cycle). The numbering of the separatrices is conducted from the selected boundary trajectory. It also determines the orientation of the boundary components. The number of the corresponding separatrix of the saddle (which belongs to the cell) is written in parentheses.

For example, for the flow in Fig. \ref{dg}, the upper saddle, and the left trajectory of the inner boundary component, the B-code looks like this: $BDS^B\{1\}[12(3)]\{23]0[\}$.

As with manifolds without boundaries, in addition to local bifurcations in typical one-parameter families of vector fields, there is a unique type of global possible bifurcation: the saddle connection. In this case, the following options are possible: 1) both saddles are internal (internal saddle connection), 2) one is internal, and the other lies on the boundary (semi-boundary saddle connection), 3) both saddles belong to the boundary (semi-boundary saddle connection).

\subsection{Saddle Connection (Internal) SC}
\begin{figure}[ht!]
\center{\includegraphics[width=0.75\linewidth]{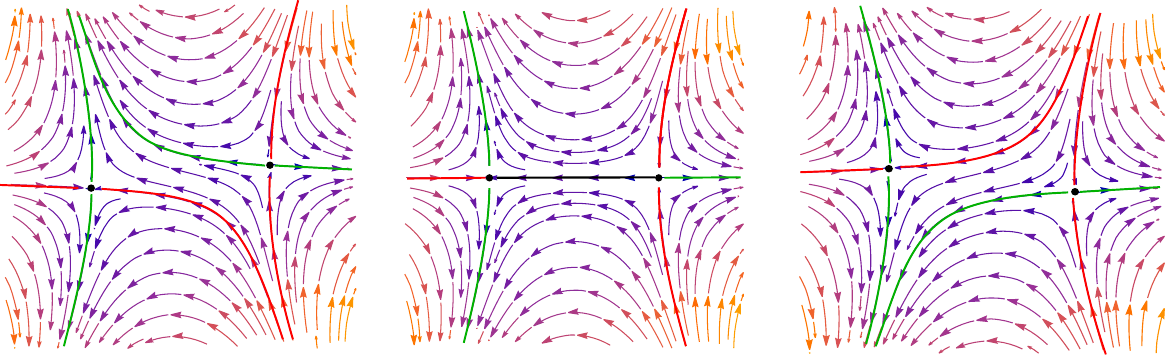}
}
\caption{saddle connection
}
\label{bifsc}
\end{figure}

The bifurcation of the saddle connection is depicted in Fig. \ref{bifsc}. It can be described by the equation $V(x,y,a)=\{x^2-y^2-1,-2xy+a\}$.
Let's describe one of the possible situations when it appears in typical one-parameter families of gradient fields. Let $p,q$ be saddle critical points of the function $f$, with $f(p)<f(q)$, and the interval $[\frac{f(p)+f(q)}{2}-\varepsilon, \frac{f(p)+f(q)}{2}+\varepsilon]$ does not contain critical values of $f$. Let $L$ be the component of the level line $f^{-1}(\frac{f(p)+f(q)}{2})$, $u$ be the separatrix for $p$, and $v$ be the separatrix for $q$, which intersect $L$.

According to Morse theory, $f^{-1}([\frac{f(p)+f(q)}{2}-\varepsilon, \frac{f(p)+f(q)}{2}+\varepsilon])$ is homeomorphic to the cylinder $S^1 \times [0,1]$, and in the coordinates $(s,t)$ (where $s$ is the polar angle on the circle), the trajectories are given by the segments $s=\text{const}$. Let us consider the twisting of the cylinder defined by the formula $(s,t) \to (s+a,t)$ (where $a$ is the deformation parameter). In this case, the trajectories will transition into helical lines. With a continuous change of the parameter $a$, the intersection point of trajectory $u$ with the circle $\frac{f(p)+f(q)}{2}+\varepsilon$ will continuously rotate around it. The intersection point of trajectory $v$ with this circle remains unchanged. By the intermediate value theorem, there is a moment in time $a$ at which these points coincide. It is precisely at this moment that a bifurcation of the saddle connection occurs.

Note that the construction described above can be realized by changing the Riemannian metric (altering the angles considered perpendicular), without changing the function itself.

In Morse theory, the bifurcation of the saddle connection corresponds to the operation of adding handles (sliding one handle over another handle of the same index).

In the diagrams, the bifurcation of the saddle connection is represented by an edge that is glued perpendicularly to the middle of another tetrahedron (in the shape of the letter T). The glued edge (1-cell), which contains the separatrix connection, is referred to as the lower edge, while the other two edges in the letter T are called the lateral edges. Thus, a flow of codimension 1  corresponds to a distinguishing graph that contains a T-vertex and the lower edge.

When writing the SN-flow code, we   start with the lower edge, which we assign the number 1. The upper edges are assigned the number 2. All other saddles and corresponding stable separatrices are numbered starting from 3 from the lower end T, following the same principle as for the Morse flow.

For example, for flow 6) in Fig. 12, the code is: $SC\{12]0[\}[2]$.

\subsection{Half-boundary saddle connection HSC}
The half-boundary saddle connection can be described by the same formula as $SC$, with a restriction on the area where it is considered. This area is formed by cutting around the separatrix connection along a stable or unstable manifold of one of the saddles that does not contain the separatrix connection. If the separatrix connection enters the boundary trajectory, we have a bifurcation $HSC_+$, and if it exits the boundary trajectory, then $HSC_-$. We obtain $hSC_+$ if we cut the bifurcation in Fig. \ref{bifsc} along the unstable manifold of the left saddle (and take the right area), while $HSC_-$ is obtained from the left area of the stable manifold of the right saddle of the bifurcation in Fig. \ref{bifsc}.

For the $HSC_+$ flow, the T-subgraph is formed by the separatrix connection and the stable manifold of the inner saddle. The base T (the lower end of the lower edge) lies on the boundary. In the code, it corresponds to 1. The inner saddle of the separatrix connection corresponds to 2.

For example, for flow 8) in Fig. 12, we have the following code: $HSC\{10[2]0\}[2]$.

\subsection{Boundary saddle connection BSC}

The boundary saddle connection $BSC$ can be described by the same formula as $SC$, with a restriction on the area where it is considered. This area is formed by cutting around the separatrix connection along the stable manifold of one saddle and the unstable manifold of another saddle, ensuring that these manifolds do not contain the separatrix connection. We obtain $BSC$ if we cut the bifurcation in Fig. \ref{bifsc} along the unstable manifold of the left saddle and the stable manifold of the right saddle, taking the central part.

In the distinguishing graph, we indicate the separatrix connection as a highlighted (black) directed edge with ends on the boundary of the surface. In the code, both ends of the separatrix connection are marked as 1, and we start writing the code from the beginning of the separatrix connection.

For example, for diagram 1) in Figure 18, the code looks like this: $BSC\{12\}\{1\}\{2]0[\}$.

\section{Topological Classification of Codimensional One Flows}
The structure of gradient flows of codimension 1 (flows at a bifurcation point) defines the structure of bifurcation. For unambiguous code writing, in addition to the starting point, it is also necessary to specify the orientation of the surface. We   list the canonical orientations defined by the flows: for $BSN$ and $BDS$, the orientation of the boundary is determined by the direction of movement at the boundary points that lie in the vicinity of the degenerate special point of the boundary. All other flows do not have canonical orientations because there exists a topological equivalence of the flow to itself, which changes the orientation (axial symmetry for flows in Figure ..). For writing codes in such cases, we use the right orientation of subsets of the plane. The code written for the opposite orientation is referred to as symmetric to the given code. If we apply axial symmetry to a compact area of the plane and the flow defined on it, we obtain a flow with a symmetric code (for the given default right orientations of the plane).

\begin{theorem} BSN and BDS flows are topologically equivalent if and only if their codes are equal. SN, SC, HN, HC, HSC, and BSC flows are topologically equivalent if and only if they have the same (unordered) pairs consisting of a code and its symmetric code.
\end{theorem}
\textit{Proof.} \textit{Necessity.} For BSN and BDS flows, there is a canonical orientation of the boundary that is opposite to the direction of motion along the trajectory of the boundary, which collapses to a point at bifurcation. This orientation is determined by the direction of motion along the trajectories adjacent to the special point at the moment of bifurcation. Therefore, the topological equivalence of such flows preserves the direction of motion along these trajectories, and thus the orientation of the surface. Consequently, it leads to the equality of codes. For the other types of bifurcations, topological equivalence may reverse the orientation. If it preserves the orientation, then the codes of the flows are equal, and their reversed codes are also equal. However, if the orientation is reversed, then the code of one flow is equal to the reversed code of the other.

\textit{Sufficiency.} If the codes of type A flows are equal, then there exists topological equivalence between the Morse flows before bifurcation, which then defines the topological equivalence of the flows at the moment of bifurcation (since they are obtained as flows on a factor space when one of the trajectories is collapsed to a point).

In the case of equal codes of type B, it is sufficient to prove that their A-codes are also equal. For the code $SN_+^B$, we increase all numbers except 0 by 1. We remove the sublist in parentheses (along with the parentheses), and in its place, we insert 1. We start $SN_+^A$ with $[1$, then write the sublist that was in parentheses with the increased numbers, close the square bracket $]$, and then write the modified list $SN_+^B$. Thus, if the codes $SN_+^B$ are equal, then the restored $SN_+^A$ codes are also equal, and the flows are topologically equivalent. For $BSN+$ bifurcation, the restoration of the A-code is similar, with the note that the $BSN+A$ code starts with $\{ [$.

$B$-codes $HS$ of flows $HN$ are essentially Morse flow codes with a fixed special point. Therefore, the topological equivalence of flows directly follows from theorems \ref{t1} and \ref{t2}.

The $BDS^B$-code transforms into a $BDS^A$-code in the following way: If a number $k$ is in parentheses, we replace $k$ with 1 and add 1 to all numbers less than $k$. We remove the round brackets and the element within them, and we add 1 as the first element (after the curly brace) to the first list, cyclically rearranging the elements in the lists so that they start with the smallest element and ordering the lists by the first element (if two lists have the same first non-zero number, their mutual order is preserved). For example, the code $BDS^B\{1\}[12(3)]\{23]0[\}$ transforms into $BDS^A\{12\}\{1]0[3\}[23]$. Thus, if the $BDS^B$-codes are equal, then the $BDS^A$-codes are also equal, and the flows are topologically equivalent.

Flow codes with a separator, just like Morse flow codes, consist of elements corresponding to the elements of the distinguishing graph, just as for Morse flows. The only difference is the presence of a separator that starts at saddle 2 and ends at saddle 1 for the codes $SC$ and $HSC_+$. Theorem, just as for Morse flows, implies that the equality of these codes leads to the equivalence of distinguishing graphs and the topological equivalence of flows.
From the construction of codes and from theorem 3, it follows

\begin{corol} Each of the codes of 1-flows has the same properties as the Morse flow codes, with the exception (in addition) of the following properties:

1) starts with one of the following fragments:

$$
\begin{array}{c}
	SN_+^A[1, \  SN_-^A[1, \ BSN_+^A \{ [1, \ BSN_-^A \{ [1
	, \ BSN_+^B\{1, \  BSN_-^B\{1, \ HS_+^A\{1,\ BSN_-^B\{1, \\ \ HS_+^B\{1, \ HS_-^B\{1, \ HN_+^A\{1,\  HN_-^A\{1,\ HN_+^B\{[1,\ HN_-^BA\{[1,\ BDS^A\{1, \ BDS^B\{1,
\end{array}
$$

2)
$BSN_+^B$ and $BSN_-^B$ contain a pair of round brackets inside square brackets.

In $HS_+^A$ and $HS_-^A$, the first unit is placed in square brackets after a cyclic permutation of the elements of the first list, while in $HS_+^B$ and $HS_-^B$, the first unit is not placed in square brackets.

In $HN_+^A$ and $HN_-^A$, the second unit is placed in the outer square brackets.

\end{corol}
\textit{Proof.}
From Theorem 4 and the definitions of codes of co-dimension 1, it follows

\begin{corol} If a sequence satisfies the properties mentioned above, then there exists a co-dimension 1 flow for which it is a code.
\end{corol}

Thus, the following types of gradient bifurcations are possible on spheres with holes:

SN -- internal saddle-node;

SC -- internal saddle connection;

BSN -- boundary saddle-node;

BDS -- boundary double saddle;

HN -- semi-boundary saddle-node (node);

HS -- semi-boundary saddle-node (saddle);

HSC -- semi-boundary saddle connection;

BSC -- connection of saddles at the boundary.

In the case of saddle-node bifurcations, such a bifurcation is defined by a separatrix diagram to the bifurcation, which highlights the trajectory (separatrix) between the saddle and the node, which collapses to a point. To define the bifurcation of a saddle connection, it is sufficient to have a separatrix diagram at the moment of bifurcation.

\section{Structure of Typical Flows and Bifurcations on the 2-Disk}

To find all possible Morse flow structures with no more than 6 singular points, we   use the following notations:

$N$ -- the number of all singular points,

$A$ -- the number of internal sources,

$B$ -- the number of boundary sources,

$S$ -- the number of internal saddles,

$T$ -- the number of boundary saddles,

$Y$ -- the number of internal sinks,

$Z$ -- the number of boundary sinks.

Note that all these numbers are non-negative and $A+B+S+T+Y+Z=N$. From the conditions of Morse flow, we also have:
$$A+B>0, \ \ \ \ Y+Z>0, \ \ \ \ B+T+Z>1, \ \ \ \ \text{B+T+Z is even,}$$
and from the Poincar\'e-Hopf theorem, it follows that
$$2A-2S+2Y+B-T+Z=\chi(DF),$$
where $\chi(DF)$ is the Euler characteristic of the double surface. In particular, for the 2-disk, $\chi(DF)=2$.

If $T_+$ is the number of boundary saddles of type $HS_+$ (with a stable separatrix), and $T_1$ is the number of boundary saddles of type $HS_+$ (with an unstable separatrix), then the following equalities hold:

$$T_+ + T_-=T, \ \ \ \ B+T_+=Z+T_-.$$

From these equalities and inequalities, we find all possible combinations for the numbers of singular points of each type for $2 \leq N \leq 6$. They are presented in the table.

\begin{center}
\begin{tabular}{|c|c|c|c|c|c|c|c|c|c|}
\hline
nn& $N$& $A$ & $B$ & $S$ & $T_+$ & $T_-$& $Y$&$Z$ & Fig.\ref{md2} \\
\hline
 1& 2&0 &1 &0 &0 &0 &0 &1 &\textbf{1}\\
\hline
2& 3& 1 & 0 & 0&1 &0 &0 &1 & 2\\
\hline
3& 4& 1 &0 &0 & 1& 1&1 & 0 & \textbf{3}\\
\hline
4 & 4 & 1 & 1 & 1 & 0 & 0 & 0 & 1& 4 \\
\hline
5 & 4 & 0 & 2 & 0& 0&1 &0 & 1& 5 \\
\hline
6 & 5 & 2 & 0 & 1& 1 & 0&0 &1 & 6 \\
\hline
7 & 5 & 1 & 1 & 1& 0& 1& 1& 0 & 7 \\
\hline
8& 5& 1 & 0 & 0& 2& 0& 0 & 2& 8 \\
\hline
9& 5 & 1 & 1 & 0 & 1& 1 & 0& 1& 9 \\
\hline
10 & 5 & 0 & 2 &1 & 0& 0& 0& 2& \textbf{10} \\
\hline
11 & 6 & 2 & 0 & 2&1 & 1& 1& 0 & 11 \\
\hline
12 &6 &2 & 1 & 2 & 0& 0& 0& 1& 12, 14  \\
\hline
13 &6 & 1 & 1 & 2& 0&0 &1 & 1& 13,\textbf{15} \\
\hline
14 &6 & 2 & 0 & 0 & 2& 1& 0& 1& 16 \\
\hline
15 &6 & 1 & 1 & 0& 1& 2& 1& 0& 17 \\
\hline
16 & 6 & 1 & 1 & 1& 1&0 & 0 & 2& 18, 20  \\
\hline
17 &6 & 1 & 2 & 1& 0& 1& 0& 1 & 19 \\
\hline
18 &6 & 0 & 2 & 0& 1 & 1& 0& 2& \textbf{21} \\
\hline
19 &6 & 0 & 1 & 0& 2& 0& 0& 3& 22 \\
\hline
\end{tabular}
\end{center}

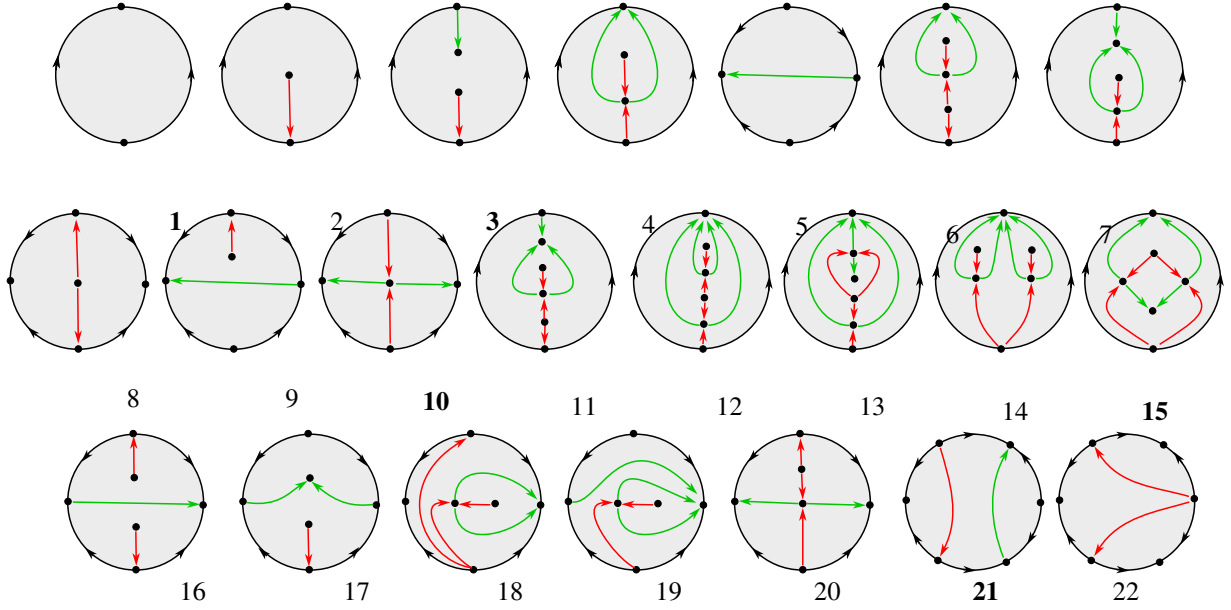
\begin{figure}[ht!]
\center{\includegraphics[width=0.98\linewidth]{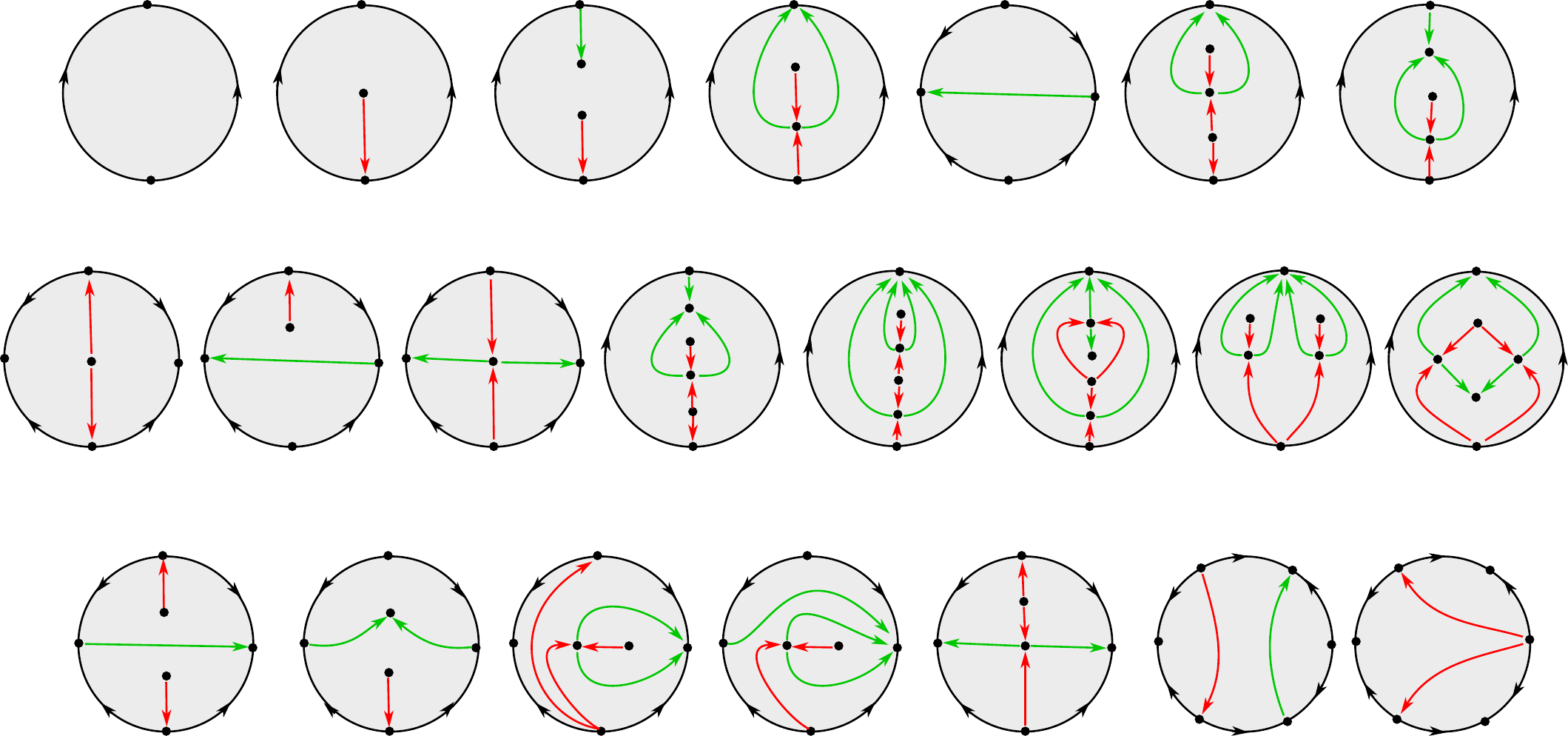}
}
\label{D35}
\put (-404,130) {\textbf{1}}
\put (-343,130) {2}
\put (-284,130) {\textbf{3}}
\put (-225,129) {4}
\put (-167,128) {5}
\put (-110,125) {6}
\put (-52,125) {7}
\put (-420,63) {8}
\put (-360,63) {9}
\put (-308,62) {\textbf{10}}
\put (-252,60) {11}
\put (-197,60) {12}
\put (-143,60) {13}
\put (-89,58) {14}
\put (-36,58) {\textbf{15}}
\put (-400,-10) {16}
\put (-338,-10) {17}
\put (-280,-10) {18}
\put (-220,-10) {19}
\put (-160,-10) {20}
\put (-100,-10) {\textbf{21}}
\put (-47,-10) {22}
\caption{Morse-flow on $D^2$ with at least 6 singular points}
\label{md2}
\end{figure}

Figure \ref{md2} shows all possible (up to homeomorphism) separatrix diagrams of Morse flows on a 2-disk with no more than 6 singular points. We   demonstrate how they are obtained using two combinations as examples.

In combination 6, there are three internal special points: two sources $v_1, v_2$ and a saddle $v_3$, as well as two boundary points: a sink $v_4$ and a saddle $v_5$ with a stable separatrix. Since there is one sink in the flow, the separatrices emanating from $v_3$ terminate at $v_4$. Therefore, the closure of their union is the boundary of a bounded region $U$. The separatrices entering $v_3$ originate from $v_1$ and $v_2$. For the sake of clarity, let $v_1 \in U$, $v_2 \in D^2 \setminus U$. Since $v_5 \in D^2 \setminus U$, the separatrix entering $v_5$ starts at $v_2$. Thus, we have uniquely determined the relative positions of all special points and separatrices. Therefore, there exists a unique separatrix diagram with this arrangement, as depicted in Figure \ref{md2}.6.

In combination 13, we have four internal points: a source $v_1$, a sink $v_2$, and two saddles $v_3, v_4$, as well as two boundary points: a source $v_5$ and a sink $v_6$. If no separatrices emanate from $v_1$, then the distinguishing graph is disconnected, which is impossible. It is also not possible for four separatrices to emanate from $v_1$, as this would cause the closure of the stable manifolds of the saddle points to divide $D^2$ into three regions, each of which must contain a sink, contradicting the fact that the flow has only two sinks. Thus, three options remain: 1) one separatrix $\gamma$ emanates from $v_1$, 2) two separatrices $\gamma_1, \gamma_2$, 3) three separatrices $\gamma_1, \gamma_2, \gamma_3$.

In the first case, let's assume for clarity that $\omega(\gamma) = v_3$. Then the remaining separatrices start at $v_5$. Let $U$ be a bounded region with boundary $\partial U = \overline{W^s(v_4)}$. Then $v_1, v_2, v_3 \in U$. Therefore, there exists a unique flow in the previous case.

In the second variant of the distinguishing graph's connectivity, it follows that $\gamma_1$ and $\gamma_2$ end at different saddles. This means that there exists a set $U\subset D^2$ such that $\partial U =\overline{W^s(v_3) \cup W^s(v_4)}$, with $v_2\in U$. Therefore, we have a unique arrangement of all singular points and separatrices, which indicates the existence of a unique flow in this variant (see Fig. \ref{md2}.15).

In the third variant, two of the three separatrices (let's say $\gamma_1$ and $\gamma_2$) end at the same saddle (let's denote it as $v_3$). Then there exists a set $U\subset D^2$ such that $\partial U =\overline{W^s(v_3)}$, with $v_2\in U$. Thus, we have a unique arrangement of all singular points and separatrices in this variant (as depicted in Fig. \ref{md2}.13). Since, in combination, the number of sources equals the number of sinks (both for the interior and the boundary), the resulting diagrams should be checked for whether they are inverses of each other. As a result, we find that the first and third variants are mutually inverse, while the second is self-inverse. Therefore, we keep two diagrams in the list (the second and third variants).

In the list of diagrams shown in Fig. \ref{md2}, the diagrams that are self-inverse are 1, 3, 10, 15, and 21. For the other diagrams, we have depicted only one of the two diagrams.

\begin{theorem}
In the two-dimensional disk, there exist the following possible structures of typical one-parameter gradient bifurcations:
\begin{itemize}
\item
with two singular points at the moment of bifurcation - two HN bifurcations;
\item
with three singular points at the moment of bifurcation: one $SN_+$, one $BSN_+$, one $HS_+$, one $HN_+$ and none of $BDS$, $SC$, $HSC_+$, $BSC$;
\item
with four singular points at the moment of bifurcation: three $SN_+$, one $BSN_+$, two $HS_+$, three $HN_+$ and none of $BDS$, $SC$, $HSC_+$, $BSC$;
\item
with five singular points at the moment of bifurcation: 12 $SN_+$, 10 $BSN_+$, 3 $HS_+$, 4 $HN_+$, 5 $BDS$, none of $SC$, 2 $HSC_+$ and none of $BSC$;
\item
with six singular points at the moment of bifurcation: 25 $SN_+$, 11 $BSN_+$, 11 $HS_+$, 22 $HN_+$, 22 $BDS$, 7 $SC$, 3 $HSC_+$ and one $BSC$.
\end{itemize}
\end{theorem}

\textit{Proof.}
First, we   show how bifurcation numbers without separatrix connections were found. If the number $n$ of special points at the moment of bifurcation is less than 6, we   use Morse flow diagrams with $n+1$ special points, where we need to indicate a separatrix or a boundary trajectory that collapses to a point at the moment of bifurcation. For example, to define the bifurcations $BSN_+$ on diagrams that are self-inverse (1, 3, 10, 15, 21), we need to indicate one red separatrix, both ends of which are internal points. On the remaining diagrams, we indicate any separatrix with internal ends. In all cases, after collapsing the specified separatrix to a point, the remaining separatrices must not form loops or oriented closed contours (as this is not possible for gradient flows). Thus, the number of possible bifurcations equals the number of separatrices that we can identify in this way. We must also consider possible symmetries of the diagram: if there is a symmetry of the diagram that transforms one identified separatrix into another, we count it only once (for example, in diagrams 14, 15). Such red separatrices are found in diagrams 4 (one), 6 (two), 7 (one), 11 (two), 12 (three), 13 (one), 15 (one), 18 (one), 18 (one), 20 (one), while green ones are in 13 (one), 14 (one). Therefore, in diagrams 3 to 5, there is only one separatrix, in 6 to 10 there are three separatrices, and in 11 to 22 there are 12 separatrices.

Using the same scheme, we find the remaining bifurcations that are not separatrix connections.

To find bifurcations with separatrix connections, one can either 1) transform the Morse flow diagram so that it has a separatrix connection, or 2) find all diagrams with separatrix connections using the same methods we used to find the Morse flow diagrams.

In the first case, to reconstruct the Morse flow diagram, there must be a cell whose boundary contains a connected path consisting of a saddle $s_1$, its stable separatrix $\gamma_1$, a source $v=\alpha(\gamma_1)$, a separatrix or boundary trajectory $\gamma_2$ $(\alpha(\gamma_2)=v, \omega(\gamma_2)=s_2)$, a saddle $s_2$, and its unstable separatrix $\gamma_3$. A saddle connection can also be obtained by reversing all motion directions along the trajectories in the described construction. In this case, we remove the separatrices $\gamma_1$ and $\gamma_3$ and replace them with a separatrix connection between $s_1$ and $s_2$.

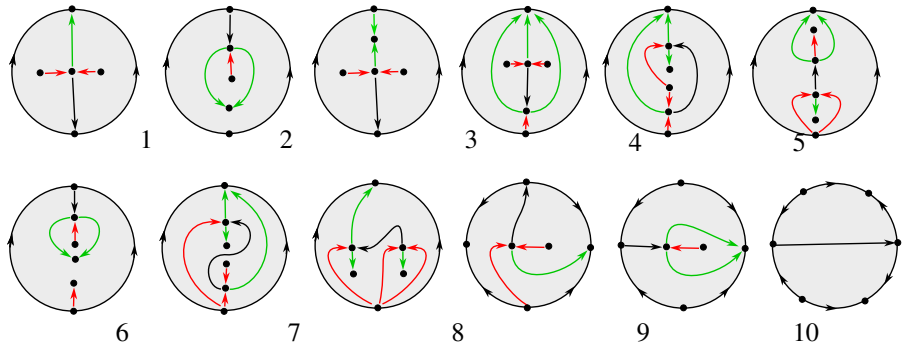
\begin{figure}[ht!]
\center{\includegraphics[width=0.72\linewidth]{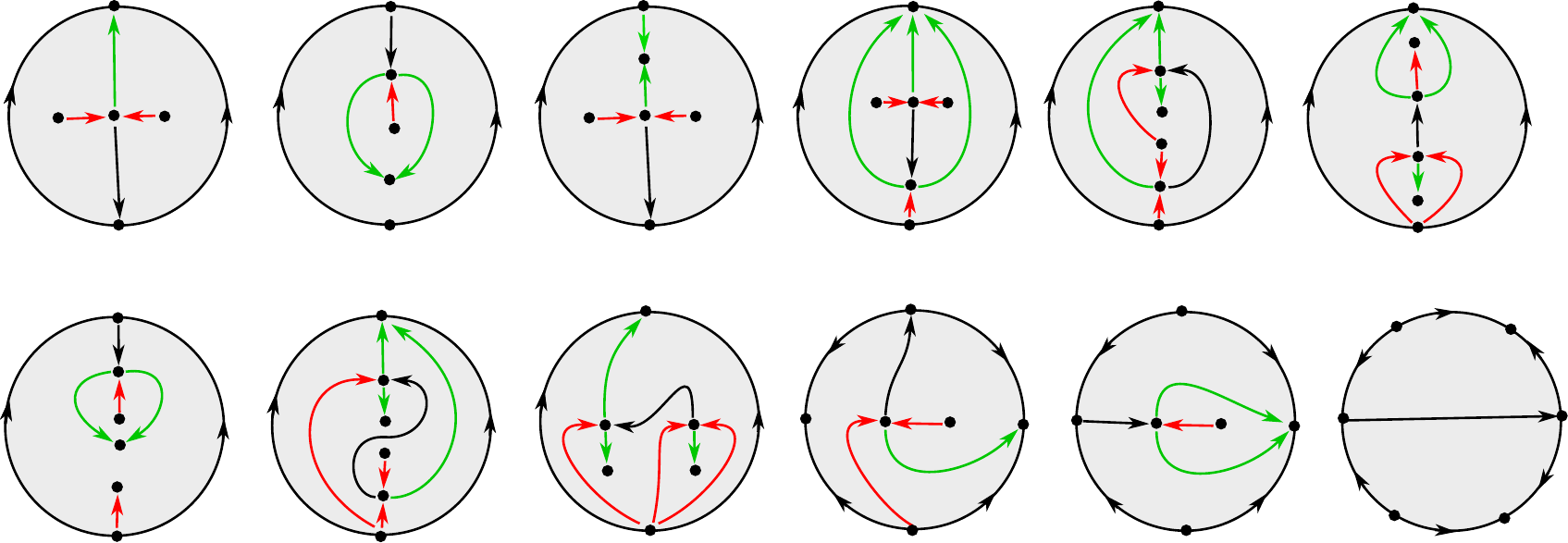}
}
\put (-288,63) {1}
\put (-235,63) {2}
\put (-165,63) {3}
\put (-103,62) {4}
\put (-41,62) {5}
\put (-297,-10) {6}
\put (-232,-10) {7}
\put (-170,-10) {8}
\put (-100,-10) {9}
\put (-41,-10) {10}
\caption{Saddle connections on $D^2$}
\label{D2sc}
\end{figure}

For example, in diagram 6 in Fig. \ref{md2}, such a path would be the saddle on the boundary (the lower vertex), the lower red separatrix, the lower internal source, the middle red separatrix, the internal saddle, and the right green separatrix. After replacing the two separatrices with a separatrix connection, we obtain diagram 1 in Fig. \ref{D2sc}.

All possible separatrix connections on a two-dimensional disk with no more than 6 singular points are depicted in Fig. \ref{D2sc}. In particular, diagrams 2, 3, and 4 are derived from diagrams 7, 11, and 12 in Fig. \ref{md2}. The diagrams of reverse flows can be obtained from these by changing all directions and colors of the separatrices (green and red swap places). With such a substitution, diagrams 6 and 10 will revert to themselves (they define one flow each). The remaining diagrams define two flows each.

\section{Structure of Typical Flows and Bifurcations on a Cylinder}

Since there must be at least two singular points on each boundary component, if the total number of such points is 4, then all of them lie on the boundary.

If the total number of singular points is 5, then one of them is internal, and there are two singular points on each boundary component (there is an even number of points on each boundary component because the sources and sinks alternate when the flow is restricted to the boundary).

For the numbers of singular points of different types, we   use the same notations as for $D^2$. Therefore, the following equalities and inequalities hold for them:
$$A+B+S+T+Y+Z=N, \ \ \ \ A+B>0, \ \ \ \ Y+Z>0, \ \ \ \ B+T+Z>3, \ \ \ \ \text{B+T+Z is even,}$$
$$2A-2S+2Y+B-T+Z=0,\ \ \ \ T_+ + T_-=T, \ \ \ \ B+T_+=Z+T_-.$$

From these equalities and inequalities, we find all possible combinations for the numbers of singular points of each type for $2 \leq N \leq$. They are presented in the table.
\begin{center}
\begin{tabular}{|c|c|c|c|c|c|c|c|c|c|}
\hline
nn& $N$& $A$ & $B$ & $S$ & $T_+$ & $T_-$& $Y$&$Z$ & Fig.\ref{SI45} \\
\hline
 1& 4&0 &1 &0 &1 &1 &0 &1 &\textbf{1,2}\\
\hline
2& 5& 1 & 0 & 0&2 &1 &0 &1 & 3\\
\hline
3& 5& 0 & 2 &1 & 0& 1&0 & 1 & 4\\
\hline
4 & 6 & 1 & 1 & 1 & 1 & 1 & 0 & 1& 5,6 \\
\hline
5 & 6 & 0 & 2 & 2& 0&0 &0 & 2& \textbf{7,8} \\
\hline
6 & 6 & 0 & 2 & 0& 1 & 2&0 &1 & 9--12 \\
\hline
\end{tabular}
\end{center}

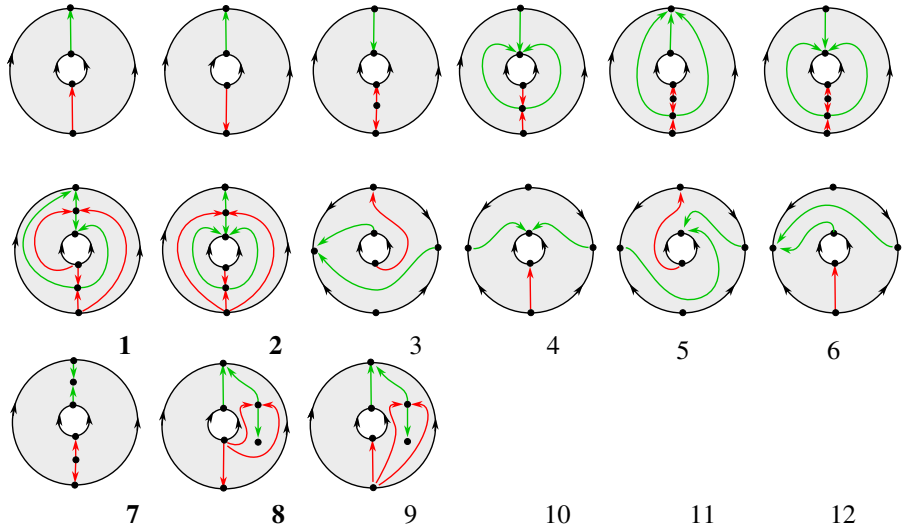
\begin{figure}[ht!]
\center{\includegraphics[width=0.72\linewidth]{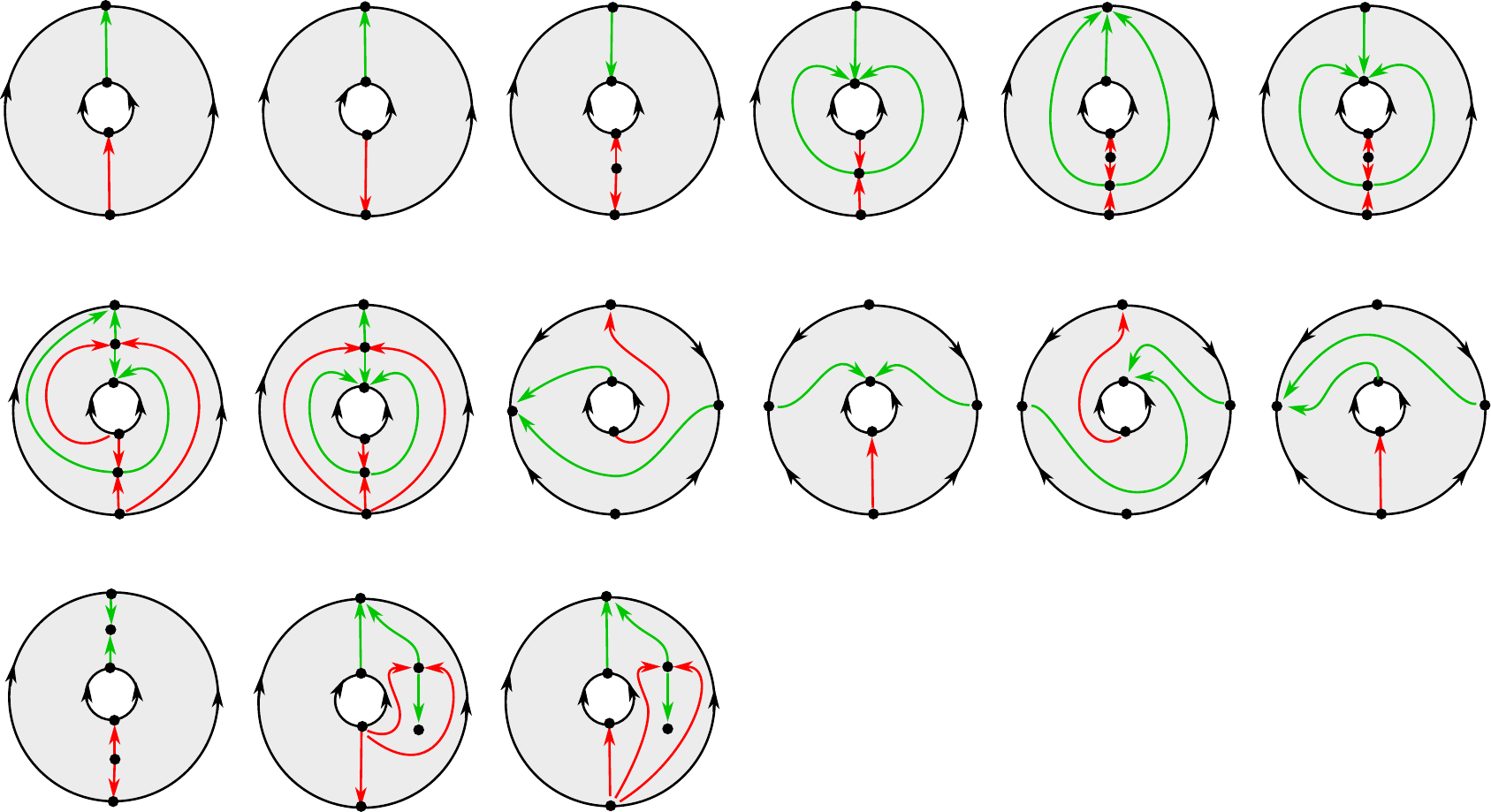}}
\put (-296,52) {\textbf{1}}
\put (-239,52) {\textbf{2}}
\put (-186,52) {3}
\put (-134,52) {4}
\put (-85,50) {5}
\put (-28,50) {6}
\put (-293,-12) {\textbf{7}}
\put (-238,-12) {\textbf{8}}
\put (-188,-12) {9}
\put (-135,-12) {10}
\put (-80,-12) {11}
\put (-27,-12) {12}
\caption{Morse flows on the cylinder }
\label{SI45}
\end{figure}

\begin{theorem}
On the cylinder, there exist two topologically non-equivalent Morse flows with four singular points (Fig.\ref{SI45}.1,2), four Morse flows with five singular points (Fig.\ref{SI45}.3,4), and fourteen Morse flows with six singular points (Fig.\ref{SI45}.5--12). Additionally, on the cylinder, the following structures of typical one-parameter gradient bifurcations saddle-node are possible:
\begin{itemize}
\item
with four singular points at the moment of bifurcation: two $HS_+$, two $HN_+$, one $BSC$;
\item
with five singular points at the moment of bifurcation: 2 $SN_+$, 7 $BSN_+$, 4 $HS_+$, 4 $HN_+$, 2 $BDS$, 0 $SC$, 1 $HSC_+$, 2 $BSC$;
\item
with six singular points at the moment of bifurcation: 16 $SN_+$, 17 $BSN_+$, 16 $HS_+$, 19 $HN_+$, 11 $BDS$, 2 $SC$, 10 $HSC_+$, 9 $BSC$.
\end{itemize}
\end{theorem}

\textit{Proof.}
In Fig. \ref{SI45}, all possible (up to homeomorphism) separatrix diagrams of Morse flows on a cylinder with no more than 6 critical points are depicted. We   demonstrate how they are obtained using the example of combination 5.

Let the flow have two internal saddles and 4 boundary points. Since the Euler characteristic of the double surface is 0, there are no saddles among the boundary points. Therefore, two of them are sources, and the other two points are sinks. For the four red separatrices, we have two options: 1) two separatrices emerge from each source (Fig. \ref{SI45}.7), 2) three separatrices emerge from one source, and one from the other (Fig. \ref{SI45}.8). The option where 4 separatrices emerge from one source is impossible, as they would then form two loops, dividing the cylinder into three regions, only two of which contain the boundary. The third internal region must have an internal sink, which contradicts the conditions for the flow.

Just as for $D^2$, to find the number of bifurcations of type $SN_+$, $BSN_+$, $BDS$ on the diagram with one more point, we find the corresponding separatrices or boundary trajectories that collapse to a point during bifurcation. For bifurcations $HS$ and $HN$ on the Morse flow diagram, we highlight the corresponding critical point on the boundary. For example, for $HS_+$, we need to highlight the saddles on the boundary with a stable separatrix for self-inverse diagrams (1,2,7,8) and both types of saddles on the boundary for the remaining diagrams. Specifically, in diagrams 1 and 2, there is one such saddle each, resulting in two bifurcations $HS_+$ with four critical points. Furthermore, in diagram 3, there are three saddles on the boundary, and in diagram 4, there is one, giving us a total of four $HS_+$ with five critical points. Finally, there are also saddles on the boundary in diagrams 5 (two), 6 (two), and 9--12 (three each). In total, we have 16 bifurcations of $HS_+$ with six critical points.
\begin{figure}[ht!]
\center{\includegraphics[width=0.95\linewidth]{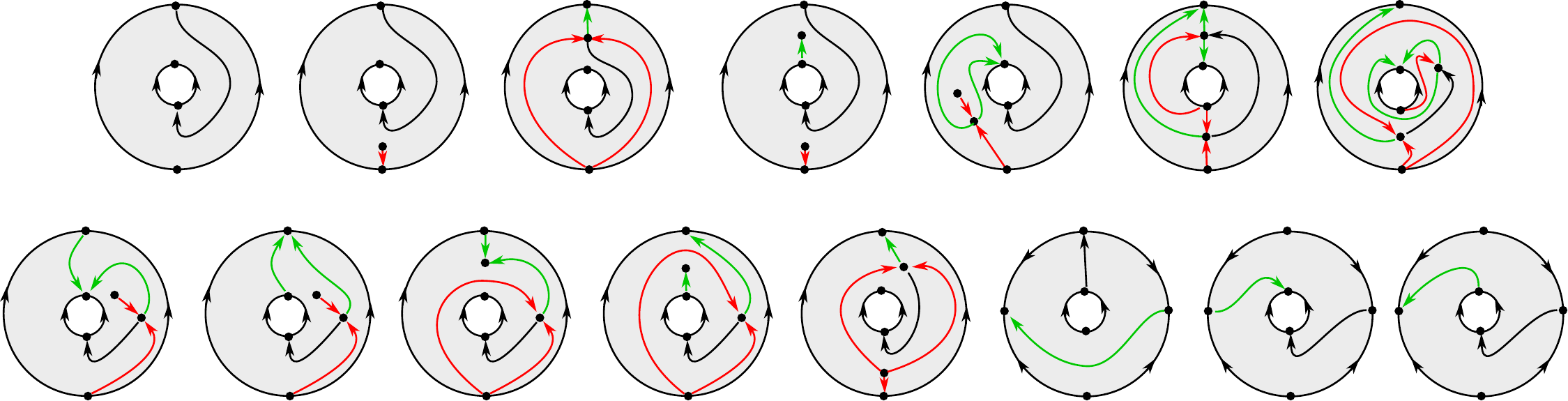}
}
\put (-381,52) {1}
\put (-327,52) {2}
\put (-273,51) {3}
\put (-212,52) {4}
\put (-157,52) {5}
\put (-105,52) {6}
\put (-53,52) {7}
\put (-403,-10) {8}
\put (-349,-10) {9}
\put (-300,-10) {10}
\put (-245,-10) {11}
\put (-192,-10) {12}
\put (-137,-10) {13}
\put (-85,-10) {14}
\put (-32,-10) {15}
\caption{Saddle connections on the cylinder}
\label{SIsc}
\end{figure}
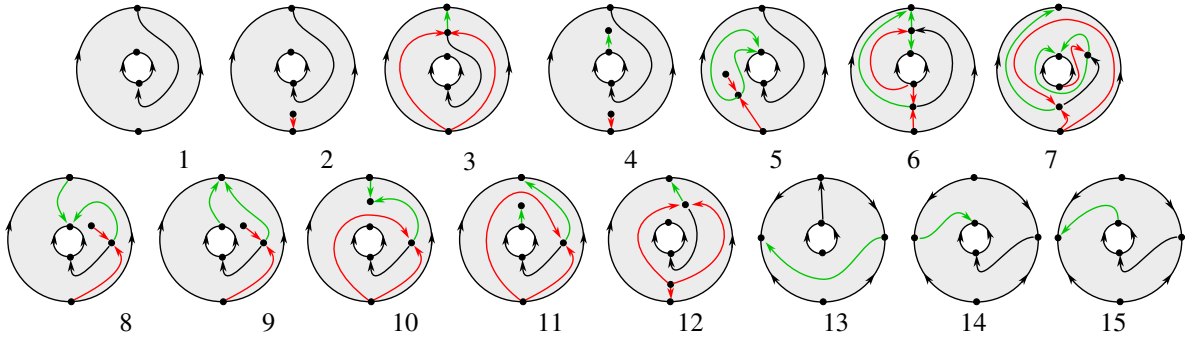

All possible saddle connections of flows with no more than six critical points on the cylinder are depicted in Fig. \ref{SIsc}. As with the cylinder, they can be obtained from Morse diagrams, where a path consisting of three trajectories is highlighted, with the first and second trajectories serving as separatrices and being replaced during the bifurcation process by a saddle connection, or by considering all flows that satisfy the same conditions as the Morse flows. For example, if a flow has 4 critical points, they lie on the boundary. Among them, there is a source, a sink, and two saddles. The saddle connection connects these two saddles. If these two saddles lie on the same boundary component, then one of the trajectories of this component together with the saddle connection forms a closed cycle, which is not possible for gradient flows. Therefore, the only remaining option is when the saddles lie on different boundary components (Fig. \ref{SIsc}.1).

Let us consider a Morse flow with 5 critical points in Fig. \ref{SI45}.3. We   replace the green and one of the red separatrices with a saddle connection. Since its ends must belong to different boundary components, we obtain the only possible diagram (Fig. \ref{SIsc}.2). By analogy, from Fig. \ref{SI45}.3 we obtain Fig. \ref{SIsc}.3.

\section{Structure of Typical Flows and Bifurcations on a Sphere with Three Holes}

We   consider the sphere with three holes as a three-punctured area on the plane (a 2-disk with two holes).

Since there must be two special points at each hole, flows with 6 special points have no other special points. Therefore, there are no internal or semi-boundary bifurcations in this case.

The doubling of the sphere with three holes results in a surface of genus 2. Besides one source and one sink, all other points are saddle points (otherwise, the Poincar\'e-Hopf theorem regarding the sum of indices would be violated).
There are only two possible structures of Morse flows in this case. In the first case, the source and sink lie on the same boundary component, while in the second case, they are on different components (see Fig. 18).

Since there are only two special points on each component of the boundary and none of the separatrices can be compressed without changing the topological type of the surface, there are no saddle-node bifurcations in this case.

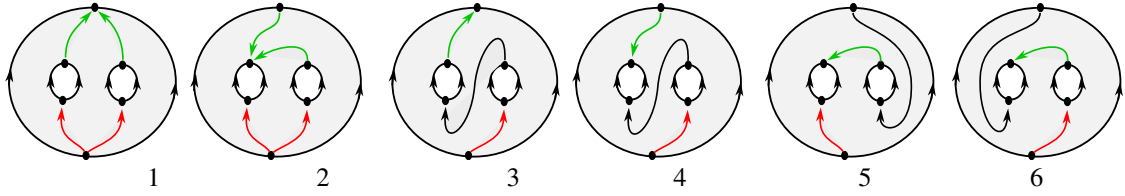
\begin{figure}[ht!]
\center{\includegraphics[width=0.9\linewidth]{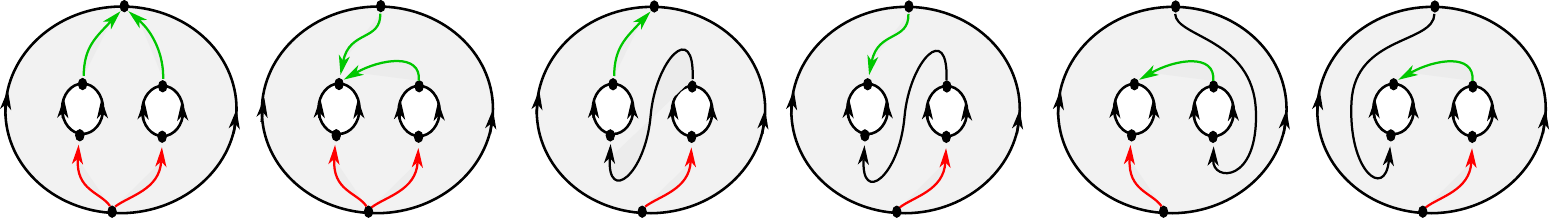}
}
\put (-369,-10) {1}
\put (-305,-10) {2}
\put (-233,-10) {3}
\put (-170,-10) {4}
\put (-100,-10) {5}
\put (-35,-10) {6}
{
\caption{ Bifurcations on the 2-sphere with tree holes }
}
\label{F03b}
\end{figure}

We obtain saddle-link bifurcations by sliding one separatrix along another. All possible 4 cases are depicted in Fig. \ref{F03b}.

Thus, the following statement is true.

\begin{theorem}
On a sphere with three holes, there are no gradient bifurcations of saddle-node type with 5 singular points at the moment of bifurcation, but there are four bifurcations of boundary saddle connections with 6 singular points.
\end{theorem}

\section*{Conclusion}

All possible structures of Morse flows and typical one-parameter bifurcations on spheres with holes containing no more than six singular points have been found. We hope that the research conducted in this work can be extended to other surfaces and with a greater number of singular points.
\begin{center}
\begin{tabular}{|c|c|c|c|c|c|c|c|c|c|}
\hline
suf,n$\setminus$bif &$SN_+$ & $BSN_+$ & $HS_+$ & $HN_+$& $BDS$ & $SC$& $HSC_+$&$BSC$ & $\sum$\\
\hline
 $D^2$, \ 3& 1&1 &1 &1 &0 &0 &0 &0 &8\\
\hline
$D^2$, \ 4& 3& 1 & 2 & 3&0 &0 &0 &0 & 18\\
\hline
$D^2$, \ 5& 12&10 &3 & 4&5 &0 &2 & 0 & 69\\
\hline
$D^2$, \ 6 & 25& 11 & 11 &22 &22 &7 &3 &1& 173 \\
\hline
$S^1\times I$, \ 4& 0&0 &2 &2 &0 &0 &0 &1 & 9\\
\hline
$S^1\times I$, \ 5&2 &7 &4 &4 &2 &0 &1 &2 & 40\\
\hline
$S^1\times I$, \ 6&16 &17 &16 &19 &11 &2 &10 &9& 168 \\
\hline
$F_{0,3}$,\ 6& 0&0 &4 &3 &0 &0 &0 &4 & 18\\
\hline
\end{tabular}
\end{center}


\section*{Acknowledgments}

We would like to express our gratitude to our collegue and students who helped us verify the accuracy of the calculations: Svitlana Bilun, Olena Myshneva and Vladislav Sinitsyn.



\bibliographystyle{unsrtnat}                                          %
\bibliography{ovtprish}
\end{document}